\documentclass[11pt]{article}
\usepackage[margin=.88in]{geometry}
\usepackage{amsmath,amssymb,amsthm,booktabs,microtype,tikz,hyperref,needspace,caption}
\usetikzlibrary{arrows.meta}
\hypersetup{colorlinks=true,linkcolor=blue,urlcolor=blue,citecolor=blue,
pdftitle={Braided Multisections and Aspherical Stein Fillings},pdfauthor={Anar Akhmedov}}
\newtheorem{theorem}{Theorem}[section]
\newtheorem{proposition}[theorem]{Proposition}
\newtheorem{lemma}[theorem]{Lemma}
\newtheorem{corollary}[theorem]{Corollary}
\theoremstyle{definition}\newtheorem{question}[theorem]{Question}
\newcommand{\Z}{\mathbb Z}
\newcommand{\Int}{\operatorname{Int}}
\newcommand{\Tor}{\operatorname{Tor}}
\newcommand{\id}{\operatorname{id}}
\numberwithin{equation}{section}
\title{Braided Multisections and Aspherical Stein Fillings}
\author{Anar Akhmedov}
\date{October 27, 2026}
\begin{document}
\maketitle
\begin{abstract}
We construct Stein fillings diffeomorphic to complements of a connected symplectic unbranched multisection and a regular sphere fiber in the ruled surface $\Sigma_g\times S^2$. For every $g\geq2$ and every multisection degree $p\geq4$, a fixed contact three-manifold admits infinitely many pairwise nonhomeomorphic aspherical fillings of this form. Their first and second integral homology groups are $\Z^{2g}\oplus\Z/p$ and $\Z^{(2g-1)(p-1)}$, respectively, and their intersection pairings vanish. We vary the braid monodromy while preserving the marked contact boundary. Unbounded torsion orders in the first homology of regular unbranched covers of degree $p2^{p-1}$ detect infinitely many isomorphism types of fundamental groups. We also give an equivariant Weinstein handle construction for the cyclic filling $W_{g,4}$. The contact boundaries are supported by spinal open books with planar pages and have positive ordinary support genus.
\end{abstract}

\section{Introduction}
We study complements of connected unbranched multisections and a regular sphere fiber in the ruled surfaces $\Sigma_g\times S^2$, with $g\geq2$.

A degree-$p$ unbranched multisection is an embedded surface that meets each sphere fiber transversely in $p$ points and projects to the base as an unbranched covering. Removing a tubular neighborhood of the multisection punctures each sphere fiber at these $p$ points. Removing a neighborhood of one sphere fiber as well leaves a bundle over the once-punctured surface $B_g=\Sigma_g\setminus\Int D^2$, with fiber the $p$-holed sphere $P_p=S^2\setminus\coprod_{j=1}^p\Int D_j^2$. Its horizontal boundary is formed by the boundary circles of the fibers, and its vertical boundary lies over $\partial B_g$. These pieces form the spine and paper of a spinal open book. The bundle admits a Stein structure whose boundary contact structure is supported by this decomposition \cite{LVW}; we recall the construction in Section~\ref{sec:cyclic}.

The multisection complements occur in the constructions of symplectic rational cohomology quadrics in \cite{AB,AY}. Their handle descriptions appear in \cite{I}; diagrams for the closed sums are being studied in \cite{II}.

We prescribe the holonomy of these surface bundles on the generators of the positive-genus base and vary the power of a Dehn twist about a separating curve while fixing the marked contact boundary. This twist acts trivially on the first homology of the planar fiber. On the mod-two homology cover of that fiber, the lift of its square acts nontrivially. Combining this cover with a cyclic cover of the base detects infinitely many isomorphism types of fundamental groups in every degree $p\geq4$. In degree four, a fiber double cover gives an eight-fold cover of the filling with explicitly computable first homology. The added twist fixes the fiber boundary collars pointwise, and the holonomy around the boundary of the base remains the identity. These two properties preserve the marked spinal open book and hence the supported boundary contact structure.

Let $r_p$ be the order-$p$ rotation of $P_p$ that cyclically permutes its boundary circles. The cyclic $P_p$-bundle over $B_g$ has monodromy $r_p$ along $a_1$ and identity along the remaining standard base generators. We denote its rounded total space by $W_{g,p}$ and its supported contact boundary by $(Y_{g,p},\xi_{g,p})$.

\begin{theorem}\label{thm:general}
For every $g\geq2$ and $p\geq4$, the contact manifold $(Y_{g,p},\xi_{g,p})$ admits aspherical Stein fillings $U_{g,p,m}$, $m\geq0$, with $U_{g,p,0}=W_{g,p}$, satisfying
\[
 H_1(U_{g,p,m};\Z)=\Z^{2g}\oplus\Z/p,\qquad
 H_2(U_{g,p,m};\Z)=\Z^{(2g-1)(p-1)}.
\]
All higher homology groups vanish, and the intersection pairing is zero. Each filling is diffeomorphic to the complement of a connected symplectic degree-$p$ unbranched multisection of genus $1+p(g-1)$ and a regular sphere fiber in a ruled surface diffeomorphic to $\Sigma_g\times S^2$.

There are positive integers $c_p,s_p$, independent of $g$ and $m$, such that for $m\geq1$ a connected regular cover of degree $p2^{p-1}$ satisfies
\[
 \bigl|\Tor H_1(\widetilde U_{g,p,m};\Z)\bigr|=c_p m^{s_p}.
\]
Consequently, for each fixed $(g,p)$, the family contains infinitely many isomorphism types of fundamental groups and an infinite pairwise nonhomeomorphic subfamily. The boundary contact structure has positive ordinary support genus.
\end{theorem}

We prove Theorem~\ref{thm:general} in Section~\ref{sec:general}. The degree-four construction gives a smaller cover and an explicit formula.

\begin{theorem}\label{thm:main}
For each $g\geq2$, the contact manifold $(Y_{g,4},\xi_{g,4})$ admits Stein fillings $V_{g,m}$, $m\geq0$, with $V_{g,0}=W_{g,4}$, such that
\[
\begin{aligned}
 H_1(V_{g,m};\Z)&=\Z^{2g}\oplus\Z/4,\\
 H_2(V_{g,m};\Z)&=\Z^{6g-3}.
\end{aligned}
\]
Every $V_{g,m}$ is aspherical and has zero intersection pairing. For $m\geq1$, it has a connected regular eight-fold cover satisfying
\[
 H_1(\widetilde V_{g,m};\Z)=\Z^{8g}\oplus(\Z/m)^2.
\]
For each fixed $g$, the family contains infinitely many isomorphism types of fundamental groups and hence an infinite pairwise nonhomeomorphic subfamily. Moreover,
\[
 H_1(Y_{g,4};\Z)=\Z^{8g-3}\oplus(\Z/4)^2,
\]
and $\xi_{g,4}$ has positive ordinary support genus.
\end{theorem}

The fillings are diffeomorphic to complements of connected symplectic multisections of genus $4g-3$ in a ruled surface diffeomorphic to $\Sigma_g\times S^2$, with one sphere fiber also removed. The contact manifolds for distinct $g$ are pairwise nondiffeomorphic, since their first Betti numbers differ.

In degree four, the general construction uses even twist powers and satisfies $U_{g,4,m}=V_{g,2m}$.

Ivan Smith \cite[Proposition~1.2]{Smith} constructed infinitely many pairwise non-isotopic connected symplectic surfaces in the class $p[\Sigma_g\times\{\mathrm{pt}\}]$ in $\Sigma_g\times S^2$ for every even degree $p\geq4$ and $g\geq1$. These surfaces have genus $1+p(g-1)$. His proof uses braid monodromy, covering spaces, and surface bundles over positive-genus bases. Chris Hays and B. Doug Park \cite{HP} studied homotopic but non-isotopic symplectic surfaces in products of Riemann surfaces. The work of Smith and of Hays and Park motivated the author's earlier work on braided multisections \cite{AB,AY}. The cyclic multisections and closed-sum gluings in those papers provide the marked ruled pieces considered here.

The author and B. Doug Park used braided symplectic surfaces to construct small exotic four-manifolds in 2007 \cite[Section~3]{APOdd}. In work carried out in 2011, the author and Weiyi Zhang explicitly used connected braided tori representing $p[T^2\times\{\mathrm{pt}\}]$ as degree-$p$ multisections of the sphere bundle $T^2\times S^2\to T^2$ in their summing-spheres construction; see \cite[Section~4, Proposition~4.1, p.~24, and Figure~4, p.~25]{AZ}. The construction joins $p$ symplectic spheres to the punctured sphere fiber across a symplectic sum along tori, producing a symplectic sphere whose self-intersection is the sum of their self-intersections. Page numbers refer to arXiv:1506.08367v2.

Akbulut and Ozbagci \cite{AO} relate section-and-fiber complements to Stein domains when the punctured fibration is positive and allowable.
Lisi, Van Horn-Morris, and Wendl \cite{LVW} developed the supported
contact and Stein structures for bordered Lefschetz fibrations.
Zhou \cite[Lemma~2.1]{Zhou} formulates the fixed-boundary filling
construction for surface bundles with compactly supported fiber
monodromy and prescribed boundary holonomy. His examples vary one
generator of a handle pair while sending the other to the identity.
Here the rotation permutes the planar fiber boundaries; the added Dehn twist fixes their collars.

Baykur and Van Horn-Morris \cite[Remark~13]{BVHM} sketch a construction of fillings of rational planar spinal open books using complements of braided surfaces in the four-ball.

The first examples of exotic Stein fillings were constructed by Akhmedov, Etnyre, Mark, and Smith \cite[Theorem~1.1]{AEMS}. They exhibited infinitely many contact three-manifolds, each admitting infinitely many homeomorphic but pairwise nondiffeomorphic simply connected Stein fillings. Akhmedov and Ozbagci \cite{AkO} extended this to arbitrary finitely presented fundamental groups: for each such group, a singularity link with its canonical contact structure admits infinitely many exotic Stein fillings with that fundamental group. They also constructed infinitely many homeomorphic but pairwise nondiffeomorphic closed symplectic four-manifolds with the prescribed fundamental group, each admitting a nonholomorphic Lefschetz fibration over $S^2$. Akbulut and Yasui \cite{AkY} constructed infinite families of exotic Stein fillings with second Betti number two. The fillings distinguished
here have different fundamental groups. Whether the contact manifolds
in Theorem~\ref{thm:main} admit homeomorphic but nondiffeomorphic Stein
fillings remains open. We plan to study this question in follow-up work.

\paragraph*{Conventions.}
All manifolds are smooth and oriented, with boundaries oriented by the outward-normal-first convention. Homology has integral coefficients unless stated otherwise, and $\mathbb F_k$ denotes the free group of rank $k$. Conjugation by a positive base generator induces the stated fiber monodromy. When comparing boundary markings, we use monodromy diffeomorphisms with specified boundary collars. Complements of surfaces mean complements of the interiors of tubular neighborhoods, with corners rounded.

\section{Cyclic multisections and Stein structures}\label{sec:cyclic}

The braided tori used by the author and Zhang in \cite[Section~4, Proposition~4.1, p.~24]{AZ} are genus-one examples of connected multisections of sphere bundles. We now give the cyclic construction over a base of genus $g\geq2$.

Let $g,p\geq2$, and let $u:\Sigma_g\to S^1$ be a smooth map representing a primitive integral cohomology class. Identify $S^1$ with the equator in $S^2$, and define
\[
 F_{g,p}=\{(x,z)\in\Sigma_g\times S^1:z^p=u(x)\}
 \subset\Sigma_g\times S^2.
\]
Projection to $\Sigma_g$ makes $F_{g,p}$ an unbranched degree-$p$ covering. Primitivity of $[u]$ makes its cyclic monodromy transitive, so $F_{g,p}$ is connected. The covering formula gives genus $G=1+p(g-1)$. A vector field normal to the equator gives a nowhere-zero normal section along $F_{g,p}$, so its self-intersection is zero. Choose a regular sphere fiber $S=\{x_0\}\times S^2$ and define
\[
 W_{g,p}=(\Sigma_g\times S^2)\setminus\Int\nu(F_{g,p}\cup S),
 \qquad Y_{g,p}=\partial W_{g,p}.
\]
Here $\nu(F_{g,p}\cup S)$ is the plumbed tubular neighborhood of the two surfaces. The tubular disks around the multisection are chosen fiberwise, and the neighborhood of $S$ is the inverse image of a small base disk. Before rounding, this gives a bundle with corners, which we round using the product collars. Write $B_g=\Sigma_g\setminus\Int D^2$ and $P_p=S^2\setminus\coprod_{j=1}^{p}\Int D_j^2$. The bundle model in \cite{I,AY} restricts to
\[
 P_p\longrightarrow W_{g,p}\longrightarrow B_g.
\]
To identify the bundle, we choose a real lift $U$ of the argument of $u$ on each simply connected base chart. Rotation through $U/p$ takes a fixed configuration of $p$ disks to tubular disks around the roots of $z^p=u(x)$. Changing $U$ by $2\pi k$ changes this trivialization by $r_p^k$. Thus the transitions factor through a character $\chi:\pi_1(B_g)\to\Z/p$. Primitivity of $[u]$ makes $\chi$ surjective and allows a standard surface marking with $u_*(a_1)=1$ and $u_*$ zero on the remaining generators. Choosing the cyclic generator with the corresponding sign gives monodromy $r_p$ along $a_1$ and identity along every other generator.

The degree-$p$ cover associated to $\ker\chi$ is
\[
 \widehat W_{g,p}=\widehat B_g\times P_p.
\]
The boundary word of $B_g$ is a product of commutators and hence has trivial image under $\chi$. Its inverse image therefore has $p$ boundary components, each mapping with degree one. If $\widehat g$ is the genus of $\widehat B_g$, then
\[
 2-2\widehat g-p=p(1-2g),
\]
so $\widehat g=1+p(g-1)=G$.

Akbulut and Ozbagci \cite[Theorem 1]{AO} consider the complement of a section and a regular fiber in a Lefschetz fibration over $S^2$. When the resulting punctured fibration over $D^2$ is positive and allowable, its total space is Stein. Here we remove a connected multisection and a sphere fiber from a ruled surface. The remaining base is $B_g$, and the boundary carries a spinal open book.

For a surface bundle over a bordered surface, the horizontal face is the bundle of fiber boundaries and the vertical face lies over the base boundary. The horizontal face is called the spine; its circle fibers are the boundary circles of the planar pages. The vertical face is the paper, fibered over $S^1$. These fibrations, together with their identifications on the common tori, form the spinal open book.

The spine is a circle bundle over the connected degree-$p$ cover $\widehat B_g\to B_g$. This cover has genus $G$ and $p$ boundary components. Since $H^2(\widehat B_g;\Z)=0$, this oriented circle bundle is trivial. We retain the inherited normal framing to specify its identifications along the boundary tori. The paper is $P_p\times S^1$, since the monodromy around $\partial B_g$ is the identity. The $p$ boundary tori are joined using the inherited product markings.

A contact form supports this decomposition if its Reeb field is positively tangent to the spine circles and its differential is positive on the paper pages, with compatible boundary orientations. The existence and uniqueness results for supported contact structures and allowable bordered Lefschetz fibrations in \cite[Theorems~1.7 and~1.24]{LVW} apply here: the fiber $P_p$ has boundary and there are no critical points. The proof below includes the compatible smoothing and passage to a Stein structure.

\begin{proposition}\label{prop:stein}
The rounded complement $W_{g,p}$ admits a Stein structure. Its boundary contact structure $\xi_{g,p}$ is supported by the spinal open book induced by the bundle, with planar page $P_p$.
\end{proposition}
\begin{proof}
We equip $B_g$ with a complex structure extending across its boundary. The cover $\widehat B_g$ inherits a complex structure with holomorphic deck transformations. Realize $P_p$ as a bordered domain in $\mathbb{CP}^1$ invariant under $r_p(z)=e^{2\pi i/p}z$. The bundle is the quotient
\[
 E=(\widehat B_g\times P_p)/(\Z/p),
\]
where the group acts by a deck transformation on the first factor and by $r_p$ or its inverse on the second, according to the transition convention. The action is free because its action on the base cover is free. Hence $E$ has an integrable complex structure.

On each bordered factor we choose a smooth strictly subharmonic function, negative in the interior and zero on the boundary. For example, solve $\Delta f=1$ with zero Dirichlet boundary value for a conformal metric. The maximum principle gives negativity and the Hopf lemma gives positive outward derivative. Averaging over the finite group preserves these properties and strict subharmonicity. We denote these invariant functions by $f_B$ and $f_P$.

On the quotient with corners we use the additive potential
\[
 f_0=f_B+f_P,\qquad
 \lambda_0=-df_0\circ J=\lambda_B+\lambda_P,
 \quad \lambda_B=-df_B\circ j_B,\quad\lambda_P=-df_P\circ j_P.
\]
These forms are invariant and descend to the quotient. The form $d\lambda_0$ is the product K\"ahler form, and its Liouville vector field is the sum of the two factor Liouville vector fields. It points outward on both boundary faces.

On the horizontal face $\widehat B_g\times\partial P_p$, we have
\[
 d\lambda_0|_{T(\widehat B_g\times\partial P_p)}
 =\operatorname{pr}_B^*(d\lambda_B).
\]
The kernel is precisely the fiber-circle direction. Since $\lambda_P$ is positive on each oriented boundary circle, the restriction of $\lambda_0$ is a contact form whose Reeb field is positively tangent to the spine fibers. On the vertical face $\partial\widehat B_g\times P_p$, the restriction of $d\lambda_0$ is the fiber area form $d\lambda_P$. It is positive on the paper pages, and $\lambda_B$ is positive in the base-boundary direction. These properties pass to the quotient and verify the spinal Giroux-form conditions.

The projection is holomorphic, and $f_0$ is constant on every connected boundary circle of each fiber. The complex tangent planes of the horizontal face are the base tangent planes; the spine Reeb flow preserves them because its coefficient depends only on the fiber-circle coordinate. Thus $(J,f_0)$ is a supported almost Stein structure in the sense of \cite[Definition 1.23]{LVW}. By \cite[Proposition~2.28]{LVW}, a compatible smoothing of the corners gives an almost Stein filling whose boundary contact structure is supported by a spinal open book with smooth overlap isotopic to the original boundary decomposition. The correspondence discussed after \cite[Definition~1.22]{LVW}, using \cite{CE}, gives a Stein structure on the rounded manifold with boundary contact structure isotopic to this one. We denote this contact structure by $\xi_{g,p}$. The uniqueness assertion of \cite[Theorem~1.7]{LVW} applies to the resulting spinal open book with smooth overlap.
\end{proof}

The complex structure is constructed on the complement. A direct strictly pseudoconvex realization of its smooth diffeomorphism type is given in Appendix~\ref{app:quotient}.

\section{Topology and contact boundaries}\label{sec:topology}

\begin{proposition}\label{prop:topology}
The complement is aspherical and satisfies
\[
\begin{aligned}
 \pi_1(W_{g,p})&\cong\mathbb F_{p-1}\rtimes_\rho\mathbb F_{2g},\\
 H_1(W_{g,p};\Z)&\cong\Z^{2g}\oplus\Z/p,\\
 H_2(W_{g,p};\Z)&\cong\Z^{(2g-1)(p-1)},\\
 H_k(W_{g,p};\Z)&=0\quad(k\geq3).
\end{aligned}
\]
Here $\rho(a_1)=(r_p)_*$ on the fiber group based at the north pole, and $\rho$ sends all remaining standard base generators to the identity. Its intersection pairing is zero, and
\[
 \chi(W_{g,p})=(2g-1)(p-2),\qquad \sigma(W_{g,p})=0.
\]
It admits a smooth handle decomposition with counts
\[
 (h_0,h_1,h_2,h_3,h_4)
 =(1,\,2g+p-1,\,2g(p-1),\,0,\,0).
\]
\end{proposition}
\begin{proof}
The rotation fixes the poles of the sphere. Either pole gives a section. Since $\pi_2(B_g)=0$, the homotopy sequence of the bundle gives
\[
 1\longrightarrow\pi_1(P_p)\longrightarrow\pi_1(W_{g,p})
 \longrightarrow\pi_1(B_g)\longrightarrow1.
\]
The section splits this sequence, giving the stated semidirect product. Both base and fiber are aspherical, so the same homotopy sequence gives $\pi_k(W_{g,p})=0$ for every $k\geq2$.

For the group presentation, we choose peripheral loops $x_1,\ldots,x_p$ in cyclic order in the fiber, based at the north pole using radial paths. They satisfy $x_1\cdots x_p=1$, and rotation sends $x_j$ to $x_{j+1}$, with indices modulo $p$. In addition to the fiber relation, the bundle group has relations
\[
 a_1x_ja_1^{-1}=x_{j+1},\qquad
 [b_1,x_j]=[a_i,x_j]=[b_i,x_j]=1\quad(i\geq2).
\]
There is no base surface relation because $B_g$ is punctured. This specifies the action in the semidirect product.

To compute integral homology, we write
\[
 A=H_1(P_p;\Z)=\Z^p/\langle e_1+\cdots+e_p\rangle.
\]
The rotation cyclically permutes the $e_j$. If the class of $v\in\Z^p$ is invariant, then
\[
 r_pv-v=k(1,\ldots,1)
\]
for an integer $k$. Summing coordinates gives $pk=0$, so $k=0$. A cyclically invariant integer vector is constant and represents zero in $A$. Thus $r_p-I$ is injective. In its cokernel all $e_j$ become a single generator $e$, and the original relation becomes $pe=0$. There are no further relations, so this cokernel is $\Z/p$.

Since $B_g$ retracts onto a bouquet of $2g$ circles, its cellular complex with coefficients in $A$ is
\[
 A^{2g}\longrightarrow A,\qquad
 (v_1,\ldots,v_{2g})\longmapsto(r_p-I)v_1.
\]
There is no degree-two chain group. Injectivity of $r_p-I$ identifies the kernel with $A^{2g-1}$, of rank $(2g-1)(p-1)$; the cokernel is $\Z/p$. The homology spectral sequence of the bundle has terms
\[
 E^2_{s,t}=H_s(B_g;H_t(P_p;\Z)).
\]
Only $s,t\in\{0,1\}$ occur, so every differential vanishes. In total degree two, this identifies $H_2(W_{g,p};\Z)$ with $H_1(B_g;A)=A^{2g-1}$. In degree one it gives
\[
 0\longrightarrow\Z/p\longrightarrow H_1(W_{g,p};\Z)
 \longrightarrow\Z^{2g}\longrightarrow0.
\]
The section splits this sequence. All higher total degrees vanish. The Euler characteristic follows either from these ranks or from
\[
 \chi(W_{g,p})=\chi(B_g)\chi(P_p)=(1-2g)(2-p).
\]

For the intersection pairing, pass to the cover $\widehat B_g\times P_p$. Both factors have zero second homology, so the K\"unneth theorem identifies its second homology with $H_1(\widehat B_g;\Z)\otimes H_1(P_p;\Z)$. Intersections of product cycles are, up to sign, products of the intersection numbers in the two factors. Since $P_p$ is planar, its first-homology pairing vanishes. Thus the cover has zero intersection pairing.

Let $\operatorname{tr}$ be transfer for this degree-$p$ cover. Counting lifts of transverse intersection points gives
\[
 Q_{\widehat W}(\operatorname{tr}x,\operatorname{tr}y)=pQ_W(x,y).
\]
The left side vanishes, and hence $Q_W=0$. In particular $\sigma(W)=0$ and the entire second homology lies in the radical of the pairing.

We use a handle decomposition of $B_g$ with one 0-handle and $2g$ 1-handles, and a decomposition of $P_p$ into one disk and $p-1$ bands. We choose the disk to be a small $r_p$-invariant neighborhood of the north pole and use the smooth handlebody conventions of \cite[Chapter~1]{AkBook}. The bundle over the base 0-handle contributes one 0-handle and $p-1$ 1-handles. Over each base 1-handle, the relative product attachment contributes one 1-handle from the fiber disk and $p-1$ 2-handles from the fiber bands. The endpoint identifications retain the specified monodromy; they need not preserve the chosen fiber decomposition. Since there is no base 2-handle, these are the required smooth handle counts.
\end{proof}

For comparison, let
\[
 M_{g,p}=(\Sigma_g\times S^2)\setminus\Int\nu F_{g,p}
\]
be the complement of the multisection alone. Over the closed base $\Sigma_g$, the local-coefficient complex has a degree-two group $A$. Tracing the boundary word $[a_1,b_1]\cdots[a_g,b_g]$ gives
\[
 \partial_2(v)=(0,(r_p-I)v,0,\ldots,0)
\]
in the coordinates $(a_1,b_1,\ldots,a_g,b_g)$, up to a common orientation sign. The degree-one differential is the same as before. It follows that
\[
 H_1(\Sigma_g;A)=\Z/p\oplus A^{2g-2}.
\]
In total degree two the bundle spectral sequence has this group and the base class $H_2(\Sigma_g;\Z)=\Z$. The pole section preserves and splits off the base class, so
\[
 H_2(M_{g,p};\Z)=\Z^{1+(2g-2)(p-1)}\oplus\Z/p.
\]
This recovers the torsion calculation in \cite{I}. A compact Stein surface has the homotopy type of a two-dimensional CW complex, whose second homology is a subgroup of its free cellular two-chain group and is therefore torsion-free. Hence $M_{g,p}$ admits no Stein structure.

Deleting the base disk removes this degree-two local-coefficient differential, and the torsion disappears from $H_2(W_{g,p};\Z)$.

The configuration $F_{g,p}\cup S$ has two square-zero components, of genera $G$ and zero, meeting positively in $p$ points. Hence $Y_{g,p}$ is the orientation reverse of the boundary of its plumbing neighborhood $N$. The plumbing graph has two vertices joined by $p$ edges. Its intersection matrix is
\[
 Q_N=\begin{pmatrix}0&p\\p&0\end{pmatrix}.
\]
The neighborhood $N$ retracts onto the union of the two surfaces. Identifying two connected surfaces at $p$ points contributes $p-1$ graph cycles, besides the $2G$ cycles on the multisection. Consequently
\[
 H_1(N;\Z)=\Z^{2G+p-1},\qquad
 H_2(N;\Z)=\Z^2,\qquad H_3(N;\Z)=0.
\]
The two surface classes generate $H_2(N;\Z)$. Because $H_1(N;\Z)$ is free, universal coefficients and Poincar\'e--Lefschetz duality identify
\[
 H_2(N,\partial N;\Z)=H^2(N;\Z)=\Z^2,
 \qquad H_1(N,\partial N;\Z)=H^3(N;\Z)=0.
\]
In the component basis and its dual, the map $H_2(N)\to H_2(N,\partial N)$ is $Q_N$. Its determinant is $-p^2$, so it is injective, and its cokernel is $(\Z/p)^2$. The pair sequence gives
\[
 0\longrightarrow\operatorname{coker}Q_N
 \longrightarrow H_1(Y_{g,p};\Z)
 \longrightarrow\Z^{2G+p-1}\longrightarrow0.
\]
The quotient is free, so the sequence splits noncanonically as a sequence of abelian groups. Hence
\begin{equation}\label{eq:boundaryH1}
 H_1(Y_{g,p};\Z)\cong\Z^{1+p(2g-1)}\oplus(\Z/p)^2.
\end{equation}

\begin{corollary}\label{cor:nonplanar}
The contact manifold $(Y_{g,p},\xi_{g,p})$ is Stein fillable but admits no supporting ordinary planar open book.
\end{corollary}
\begin{proof}
The zero intersection pairing and the homology calculation give $b_2^0(W_{g,p})=(2g-1)(p-1)>0$, where $b_2^0$ is the dimension of the radical of the absolute intersection pairing. If $\xi_{g,p}$ admitted an ordinary planar supporting open book, Etnyre \cite[Theorem 4.1]{Etnyre} would force $b_2^0=0$ for every symplectic filling, including the Stein filling just constructed. Hence the support genus is at least one.
\end{proof}

The three pieces in \cite{I} have the following invariants. The last column records smooth handle counts.
\begin{center}
\begin{tabular}{ccccc}
\toprule
$(g,p)$ & $H_1(W_{g,p};\Z)$ & $b_2(W_{g,p})$ & $\chi(W_{g,p})$ & $(h_0,h_1,h_2)$\\
\midrule
$(2,3)$ & $\Z^4\oplus\Z/3$ & $6$ & $3$ & $(1,6,8)$\\
$(3,2)$ & $\Z^6\oplus\Z/2$ & $5$ & $0$ & $(1,7,6)$\\
$(2,4)$ & $\Z^4\oplus\Z/4$ & $9$ & $6$ & $(1,7,12)$\\
\bottomrule
\end{tabular}
\end{center}
For $W_{2,3}$, the pages of the spinal open book are pairs of pants, and the base of the spine has genus four and three boundary components, and
\[
 \pi_1(W_{2,3})=\mathbb F_2\rtimes\mathbb F_4,
 \qquad H_1(Y_{2,3};\Z)=\Z^{10}\oplus(\Z/3)^2.
\]

In this case the fiber group has generators $x=x_1$ and $y=x_2$, with $x_3=(xy)^{-1}$. The rotating base generator acts by
\[
 a_1xa_1^{-1}=y,\qquad a_1ya_1^{-1}=(xy)^{-1}.
\]
The other three base generators commute with $x$ and $y$. In abelianization, the first relation gives $x=y$ and the second gives $y=-x-y$, hence $3x=0$. This recovers the summand $\Z/3$ in $H_1(W_{2,3};\Z)$. Its six second-homology generators lie in the radical of the intersection pairing. The corresponding boundary calculations for the other two examples give
\[
 H_1(Y_{3,2};\Z)=\Z^{11}\oplus(\Z/2)^2,\qquad
 H_1(Y_{2,4};\Z)=\Z^{13}\oplus(\Z/4)^2.
\]

\section{Equivariant Weinstein handles}\label{sec:handles}

We use a rotation-invariant ribbon decomposition of the fiber, with disks $D_N,D_S$ at the poles and four bands $e_0,e_1,e_2,e_3$ joining them. Their cyclic orders at the two vertex disks are opposite. The resulting ribbon surface is $P_4$, and rotation fixes the vertex disks and sends $e_j$ to $e_{j+1}$, with indices modulo four.

\begin{center}
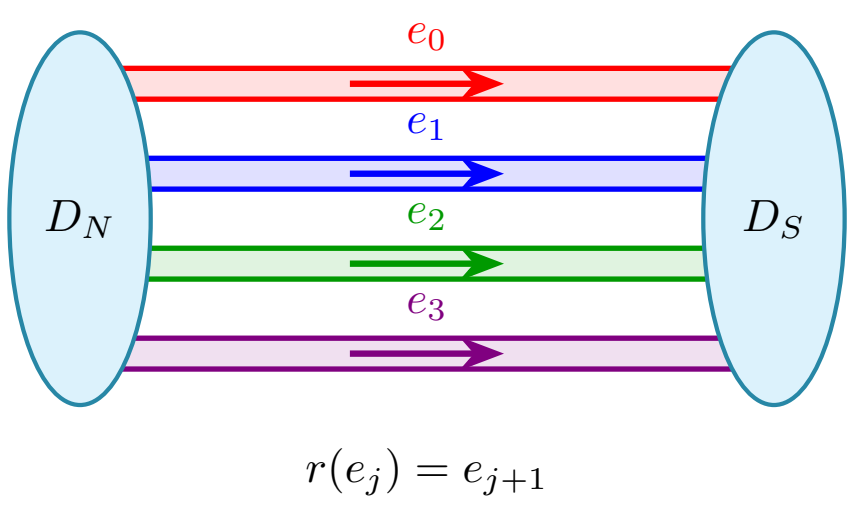

\begin{minipage}{\linewidth}
\centering
\begin{tikzpicture}[scale=.85,transform shape,x=1cm,y=1cm,>=Stealth,font=\small]
\foreach \y/\col/\j in {1.05/red/0,.35/blue/1,-.35/green!60!black/2,-1.05/violet/3}{
 \filldraw[fill=\col!12,draw=\col,line width=1pt] (-2.7,\y-.12) rectangle (2.7,\y+.12);
 \draw[\col,line width=1.2pt,->] (-.6,\y)--(.6,\y);
 \node[text=\col,above] at (0,\y+.13) {$e_{\j}$};
}
\filldraw[fill=cyan!12,draw=cyan!60!black,thick] (-2.7,0) ellipse (.55 and 1.45);
\filldraw[fill=cyan!12,draw=cyan!60!black,thick] (2.7,0) ellipse (.55 and 1.45);
\node at (-2.7,0) {$D_N$};\node at (2.7,0) {$D_S$};
\node[below] at (0,-1.65) {$r(e_j)=e_{j+1}$};
\end{tikzpicture}

\captionof{figure}{The invariant fiber ribbons, oriented from north to south.}\label{fig:fiber}
\end{minipage}
\end{center}

\begin{lemma}\label{lem:invariantfiber}
The four-holed sphere $P_4$ admits an $r$-invariant Weinstein surface
structure with two index-zero critical points at the poles and four
index-one critical points permuted by $r$. Its Lyapunov function may be
chosen to have value zero at the minima, value one at the saddles,
and constant value three on the boundary.
\end{lemma}
\begin{proof}
We start with the invariant graph consisting of the two poles and four meridional
arcs joining them. A regular neighborhood of this graph has two vertex
disks and four bands. The cyclic orders of the bands at the two disks
are opposite, so the thickening is a planar surface with four boundary
components. Its complementary disks form one free orbit under rotation.
Choosing the neighborhood and its boundary collars equivariantly
identifies this surface with $P_4$.

We use the surface Weinstein handle construction of \cite{CE}, starting with the standard radial structure on each vertex disk. In
positively oriented coordinates centered at the north and south poles,
the sphere rotation acts through angles $\pi/2$ and $-\pi/2$,
respectively. Both preserve the radial Liouville form. Choose four
equally spaced attaching intervals on each disk, following the
meridional graph, and attach four copies of the same Weinstein
one-handle. The attaching intervals have the opposite cyclic orders
at the two disks. Choose the attaching collar maps for one handle
and obtain the other three maps by rotation. The Liouville forms
match in these collars by the standard Weinstein handle gluing
construction. Rotation extends by sending each copy of the handle
to the next using the identical handle coordinates. After four steps
this extension is the identity. Thus the resulting Liouville form
and vector field are invariant under the prescribed order-four action.

An index-one local model has symplectic form
$dq\wedge dp$, Liouville form $-q\,dp-2p\,dq$, and vector field
$-q\partial_q+2p\partial_p$. The function $1-q^2+p^2$ is strictly
increasing along this vector field away from its saddle. We use this
model on each band and identical radial functions of value zero
at the two disk centers. Choose their attaching levels below one
and their outgoing levels above one. The standard handle interpolation
gives an invariant Lyapunov function with precisely these two minima
and four saddles. Extend the outgoing collars to a common level three
and round their corners equivariantly. This introduces no critical
points and makes the Liouville vector field outward along the boundary.
\end{proof}

\begin{proposition}\label{prop:equivariantW24}
The manifold $W_{2,4}$ admits a Weinstein handle decomposition with two $0$-handles, twelve $1$-handles, and sixteen $2$-handles. It has a Stein realization whose boundary contact structure is isotopic to $\xi_{2,4}$. One connecting $0/1$-handle pair can be cancelled, leaving one $0$-handle, eleven $1$-handles, and sixteen $2$-handles. The $2$-handles are attached along Legendrian knots with framing $-1$ relative to the contact framing.
\end{proposition}
\begin{proof}
We use the fiber data from Lemma~\ref{lem:invariantfiber}.
We choose a Weinstein surface structure on $B_2$ with one minimum
and four saddles. As in the surface handle construction above,
we choose its Lyapunov function to have value zero at the minimum,
value one at the saddles, and constant value three on the boundary.
We lift these data to the cyclic cover $\widehat B_2$. The lifted
structure has four minima and sixteen saddles, with the same
critical values. Its deck transformations preserve the Liouville
form and the Lyapunov function.

On $\widehat B_2\times P_4$, we use the product Weinstein structure
\[
 \lambda=\lambda_B+\lambda_P,\qquad
 Z=Z_B+Z_P,\qquad
 \phi=\phi_B+\phi_P.
\]
The function $\phi$ is strictly increasing along $Z$ away from the
product critical points, and its Hessian is the direct sum of the
factor Hessians at each critical point. Thus the indices add.
The free diagonal cyclic action preserves these data, so they
descend to the quotient. The Liouville vector field is outward on both
boundary faces. We round the corners equivariantly and adjust the
Lyapunov function only in an outward boundary collar, using the
Liouville flow, so that the rounded boundary is a regular level.
No critical point is introduced. This adjustment can be kept away
from $\{\phi\leq2+\eta\}$ for a fixed small $\eta>0$: each original
boundary face has $\phi\geq3$, whereas every product critical value
is at most two. The quotient is therefore a Weinstein domain with
the same interior critical points.

On the cover their numbers in indices zero, one, and two are
\[
 4\cdot2=8,\qquad
 16\cdot2+4\cdot4=48,\qquad
 16\cdot4=64.
\]
Every critical-point orbit has four elements. Hence the quotient
has handle counts $(2,12,16)$. All critical points of indices
zero, one, and two have values zero, one, and two, respectively.
We choose a regular level $1<c<2$. Its sublevel contains exactly the
subcritical handles. The stable manifold, for the forward Liouville
flow, of each index-two critical point is a Lagrangian core disk;
in the product model it is the product of the two one-dimensional
factor cores. Truncating this stable manifold at $\phi=c$ gives its
Legendrian attaching circle in the contact boundary of the
subcritical domain. Indeed, $Z$ is tangent to the disk and the
symplectic form vanishes on it, so $\lambda=\iota_Zd\lambda$
vanishes on its tangent bundle. The standard Weinstein attachment
has smooth framing one less than the contact framing.
Gompf's Stein handle attachment theorem therefore gives a
compatible Stein realization \cite[Proposition~2.3]{Gompf}; see also \cite[Section~8.2]{AkBook}.

For the saddle corresponding to $e_0$ and a base minimum, the two descending separatrices have constant base coordinate, one ending at each vertex minimum. The free deck action identifies no further trajectories. The north and south minima remain distinct in the quotient, so exactly one trajectory connects this saddle to the south minimum. Along this trajectory, the stable manifold of the saddle intersects the unstable manifold of the south minimum transversely: the latter is an open four-dimensional submanifold. Thus the unique connecting trajectory satisfies the transversality condition for Weinstein $0/1$-handle cancellation \cite{CE}. The higher attachments are carried along; the cancellation need not be equivariant.

On $\widehat B_2\times\partial P_4$, the restriction of $d\lambda$ is the base area form and its kernel is the fiber-circle direction. Since the fiber boundary Liouville form is positive, the Reeb field is positively tangent to the spine circles. On $\partial\widehat B_2\times P_4$, it is the fiber area form, positive on the paper pages. The Liouville form is positive on the corresponding oriented boundary directions. These facts descend, and the Liouville vector field is outward along both faces. After the supported corner rounding, the resulting contact structure is supported by the original marked spinal open book. The uniqueness theorem \cite[Theorem 1.7]{LVW} identifies it with $\xi_{2,4}$.
\end{proof}

We write $a=a_1$, $b=b_1$, $c=a_2$, and $d=b_2$. Before cancellation, the base bands give eight $1$-handles $t_N,t_S$, for $t\in\{a,b,c,d\}$. The four fiber bands $e_j$ give the remaining $1$-handles. We use $\epsilon_a=1$ and $\epsilon_b=\epsilon_c=\epsilon_d=0$. With the convention that conjugation by $a$ induces $r_*$, the product $2$-handle $K_{t,j}$ has attaching word
\begin{equation}\label{eq:W24routes}
 t_N e_j t_S^{-1}e_{j+\epsilon_t}^{-1},\qquad j=0,1,2,3.
\end{equation}

\begin{center}
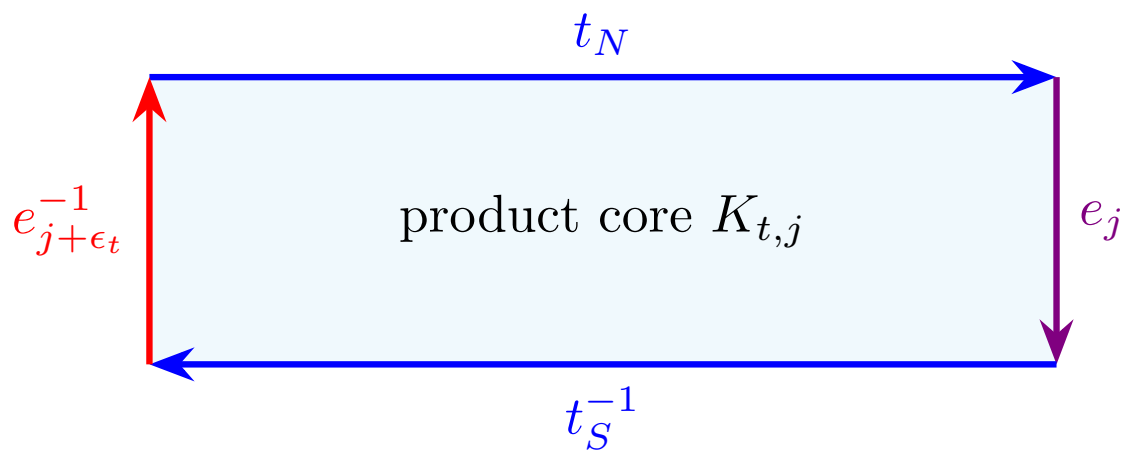

\begin{minipage}{\linewidth}
\centering
\begin{tikzpicture}[x=1cm,y=1cm,>=Stealth,font=\small]
\fill[cyan!5] (0,0) rectangle (6,1.9);
\draw[blue,very thick,->] (0,1.9)--node[above] {$t_N$}(6,1.9);
\draw[violet,very thick,->] (6,1.9)--node[right] {$e_j$}(6,0);
\draw[blue,very thick,->] (6,0)--node[below] {$t_S^{-1}$}(0,0);
\draw[red,very thick,->] (0,0)--node[left] {$e_{j+\epsilon_t}^{-1}$}(0,1.9);
\node at (3,.95) {product core $K_{t,j}$};
\end{tikzpicture}

\captionof{figure}{The attaching loop of a product $2$-handle.}\label{fig:route}
\end{minipage}
\end{center}

The product coordinates determine the ribbon embeddings and the smooth attaching framings. The coefficient $-1$ in Proposition~\ref{prop:equivariantW24} is measured relative to the contact framing.

After cancelling $e_0$ against the south vertex $0$-handle, we write
\[
 z_j=e_je_0^{-1},\quad z_0=1,\qquad
 t=t_N,\quad \bar t=e_0t_Se_0^{-1}.
\]
These band generators are distinct from the peripheral generators used earlier. The remaining eleven $1$-handles are
\[
 a,b,c,d,\bar a,\bar b,\bar c,\bar d,z_1,z_2,z_3,
\]
and the sixteen attaching words are $t z_j\bar t^{-1}z_{j+\epsilon_t}^{-1}$. For $t=b,c,d$, the $j=0$ relation gives $\bar t=t$; the remaining relations give $[t,z_j]=1$. For $a$ they give
\[
 a=z_1\bar a,\quad az_1=z_2\bar a,\quad
 az_2=z_3\bar a,\quad az_3=\bar a.
\]
Consequently $\bar a=z_1^{-1}a$ and
\[
 az_ja^{-1}=z_{j+1}z_1^{-1},\qquad j=1,2,3,\quad z_4=z_0=1.
\]
This agrees with the action of $r$ on the band generators.

With the boundary components ordered by the ribbon surface, the peripheral generators are
\[
 x_1=z_1^{-1},\qquad x_2=z_1z_2^{-1},\qquad
 x_3=z_2z_3^{-1},\qquad x_4=z_3.
\]
Their product is one. The displayed action on the $z_j$ gives $r_*(x_j)=x_{j+1}$, with indices modulo four. These are the peripheral generators used in Section~\ref{sec:variation}.

In abelianization, $\bar a=a-z_1$, $z_2=2z_1$, $z_3=3z_1$, and $4z_1=0$. Thus the handle chain complex has
\[
 H_1=\Z^4\oplus\Z/4,\qquad H_2=\Z^9,
 \qquad\chi=1-11+16=6.
\]
The map $d_2:\Z^{16}\to\Z^{11}$ has seven nonzero Smith factors, $1,1,1,1,1,1,4$. Eliminating the four barred generators leaves seven generators and twelve relations.

The same equivariant construction for $B_g$ gives handle counts
\[
 (2,4g+4,8g,0,0)\longrightarrow(1,4g+3,8g,0,0)
\]
after one connecting $0/1$-handle cancellation.

The smaller counts in Proposition~\ref{prop:topology} are smooth handle counts. Eliminating the barred generators simplifies the presentation but does not cancel Weinstein handles; such cancellations must respect the Legendrian attachments and their framings. Figures~\ref{fig:fiber} and \ref{fig:route} describe the ribbon and product attachments, rather than a simultaneous Legendrian front.

\section{Pure-braid variation and covering homology}\label{sec:variation}

We assign a power of a Dehn twist to a second base generator of the degree-four ruled piece. This preserves the marked spinal boundary.

Fix $g\geq2$. Throughout this section we abbreviate $V_{g,m}$ to $W_m$.
Let $P=P_4$, with the rotation and based peripheral system fixed in
Section~\ref{sec:handles}. We specify the twisting curve relative to this system, as in Figure~\ref{fig:twist}.
The band cores $e_0,e_2$, joined inside the two vertex disks, form an embedded
circle. Its north-based representative is $e_0e_2^{-1}$. The side containing
the first two boundary components is a pair of pants. We push this circle slightly into
that side, leaving the north basepoint outside, and denote the resulting
curve by $\gamma$. The push-off also specifies a joining arc from the
basepoint to $\gamma$. With the corresponding orientation, its based word is
\[
 e_0e_2^{-1}=z_2^{-1}=x_1x_2=:w.
\]
The band and peripheral markings therefore determine $\gamma$ and its based joining arc. We choose
its annular neighborhood disjoint from the north pole and the four boundary collars. Let $\tau$ be a Dehn twist in this annulus, with the sign fixed below.

\begin{center}
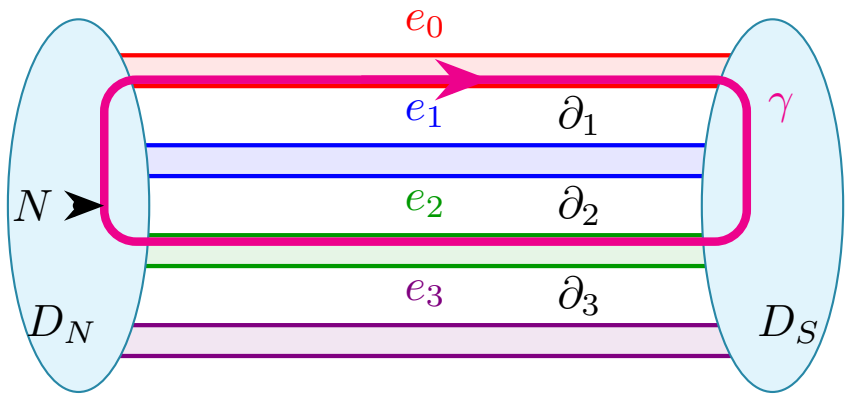

\begin{minipage}{\linewidth}\centering
\begin{tikzpicture}[scale=.85,transform shape,x=1cm,y=1cm,>=Stealth,font=\small]
\foreach \y/\col/\j in {1.05/red/0,.35/blue/1,-.35/green!60!black/2,-1.05/violet/3}{
 \filldraw[fill=\col!10,draw=\col,line width=.8pt] (-2.7,\y-.12) rectangle (2.7,\y+.12);
 \node[text=\col,above] at (0,\y+.13) {$e_{\j}$};
}
\filldraw[fill=cyan!10,draw=cyan!60!black] (-2.7,0) ellipse (.55 and 1.45);
\filldraw[fill=cyan!10,draw=cyan!60!black] (2.7,0) ellipse (.55 and 1.45);
\draw[magenta,line width=1.6pt,rounded corners=6pt]
 (-2.5,-.28) rectangle (2.5,.98);
\draw[magenta,line width=1.6pt,->] (-.5,.98)--(.5,.98);
\fill (-2.7,0) circle (.045);
\node[left] at (-2.76,0) {$N$};
\draw[black,line width=1pt,->] (-2.7,0)--(-2.5,0);
\node[text=magenta,right] at (2.53,.76) {$\gamma$};
\node at (-2.83,-.92) {$D_N$};\node at (2.83,-.92) {$D_S$};
\node at (1.2,.70) {$\partial_1$};\node at (1.2,0) {$\partial_2$};
\node at (1.2,-.70) {$\partial_3$};
\node[below] at (0,-1.58) {$w=e_0e_2^{-1}=x_1x_2$};
\end{tikzpicture}
\captionof{figure}{The standard twist curve and its joining arc from the north pole.}\label{fig:twist}
\end{minipage}
\end{center}

\begin{lemma}\label{lem:standard-twist}
The twist can be chosen to fix the basepoint and all boundary collars, with
\[
 \tau_*(x_1)=wx_1w^{-1},\qquad
 \tau_*(x_2)=wx_2w^{-1},\qquad
 \tau_*(x_3)=x_3,\qquad \tau_*(x_4)=x_4.
\]
\end{lemma}
\begin{proof}
The first two peripheral loops can be represented using the joining arc,
followed by arcs in the pair of pants bounded by $\gamma$ and the first two
boundary components. Their entry and exit arcs cross the twisting annulus once in each
direction. The other two peripheral loops can be represented outside the
annulus.

To fix the sign, use coordinates $(\theta,s)\in\mathbb R/(2\pi\mathbb Z)
\times[-1,1]$ on the annulus, with $s=-1$ on the side containing the basepoint.
Orient $\theta$ so that its based loop is $w$. Choose a smooth function
$h(s)$ equal to zero near $-1$ and one near $1$, and set
\[
 \tau(\theta,s)=(\theta+2\pi h(s),s).
\]
Extend this map by the identity. The entry arc acquires one copy of $w$ and
the exit arc acquires its inverse. This gives the two conjugation formulas.
The outside loops are fixed. Equivalently, $\tau_*(w)=w$ and $x_3$ is fixed,
so $x_4=(wx_3)^{-1}$ is fixed as well.
\end{proof}

We retain the rotation convention $r_*(x_j)=x_{j+1}$, with indices modulo
four. We use the based peripheral generators fixed above.

For $m\geq0$, we define a monodromy representation $\rho_m:\pi_1(B_g)\to\operatorname{Diff}^+(P)$ by
\[
 \rho_m(a_1)=r,\qquad \rho_m(b_1)=\mathrm{id},\qquad
 \rho_m(a_2)=\tau^m,\qquad \rho_m(b_2)=\mathrm{id}.
\]
All remaining base generators act by the identity. Let $W_m\to B_g$ be the resulting $P$-bundle, with its corners rounded. Thus $W_0=W_{g,4}$ with the markings used above.

\begin{lemma}\label{lem:boundary}
Let $B$ be a compact connected oriented surface with nonempty boundary and a fixed basepoint. Choose a parametrization of each boundary circle and a path from the basepoint to it, giving based boundary loops $\beta_i$. Let $P$ be a compact connected oriented surface with nonempty boundary, with fixed parametrized collars $C$ of its boundary circles. Suppose that actual representations
\[
 \rho_s:\pi_1(B)\longrightarrow\operatorname{Diff}^+(P)
\]
preserve $C$, have identical restrictions to $C$, and satisfy
\[
 \rho_s(\beta_i)=h_i
\]
for fixed diffeomorphisms $h_i$ of $P$. These equalities are in $\operatorname{Diff}^+(P)$, rather than only in the mapping class group. The markings retain the parametrized collars, the chosen paths to the boundary circles, and the resulting identifications of the paper with the mapping tori of the $h_i$. Then the induced marked spinal open books are isomorphic. The associated supported Stein fillings have contactomorphic boundaries under identifications in these marking classes.
\end{lemma}
\begin{proof}
Use the chosen paths to identify the fibers over the boundary basepoints with $P$. Since the collar representations agree, the horizontal boundary bundles identify by the identity in their flat collar coordinates. The paper over the $i$th boundary circle is the mapping torus of $h_i$, with its fixed fiber and base parametrizations. The equality of collar restrictions makes these identifications agree on every corner torus. They therefore give an isomorphism of the marked spinal open books.

We use the same product-collar smoothing on these boundary pieces and the same adjustment to obtain a smooth overlap as in \cite[Remark 1.6]{LVW}. Their supported contact structures are isotopic by \cite[Theorem 1.7]{LVW}. Gray stability adjusts each boundary identification through an isotopy to a contactomorphism. The surface bundles are allowable bordered Lefschetz fibrations with no critical points. Supported almost Stein structures exist by \cite[Theorem~1.24]{LVW}; the rounding and passage to a Stein structure described in Section~\ref{sec:cyclic} give the required Stein fillings.
\end{proof}

\begin{proposition}\label{prop:fixedboundary}
For every $m\geq0$, $W_m$ admits a Stein structure filling $(Y_{g,4},\xi_{g,4})$.
\end{proposition}
\begin{proof}

Since $\tau$ is the identity near the four boundary components, the horizontal boundary representations of all the bundles agree. They are the original cyclic representation on $\partial P$. The base-boundary monodromy is
\[
 \rho_m\left(\prod_{i=1}^g[a_i,b_i]\right)
 =[r,\mathrm{id}][\tau^m,\mathrm{id}]=\mathrm{id}.
\]
The paper is therefore $P\times S^1$ with the same collar identification for every $m$. On the overlaps with the horizontal boundary these trivializations agree, since the twist fixes the boundary collars pointwise. Lemma~\ref{lem:boundary} identifies the marked spinal open books for all $m$.

Lemma~\ref{lem:boundary} therefore gives a Stein structure on $W_m$ filling $(Y_{g,4},\xi_{g,4})$.

\end{proof}

\begin{proposition}\label{prop:ordinary}
Every $W_m$ is aspherical, with $\pi_1(W_m)=\mathbb F_3\rtimes_{\rho_m}\mathbb F_{2g}$, and
\[
 H_1(W_m;\Z)=\Z^{2g}\oplus\Z/4,\qquad
 H_2(W_m;\Z)=\Z^{6g-3},\qquad H_k(W_m;\Z)=0\quad(k\geq3).
\]
Its intersection pairing is zero and $\chi(W_m)=4g-2$.
\end{proposition}
\begin{proof}
Both holonomies fix the north pole, so the bundle has a section. Its fundamental group is $\mathbb F_3\rtimes_{\rho_m}\mathbb F_{2g}$, and its higher homotopy groups vanish. The twist acts trivially on $A=H_1(P;\Z)$, since the intersection pairing of a planar surface is zero. Thus the local-coefficient differential on the base graph is still
\[
 A^{2g}\longrightarrow A,\qquad(v_1,\ldots,v_{2g})\longmapsto(r_*-I)v_1.
\]
The local-coefficient calculation in Proposition~\ref{prop:topology} gives the homology and Euler characteristic.

For the intersection pairing, the fixed boundary has $b_1(Y_{g,4})=8g-3$, while $b_1(W_m)=2g$ and $b_2(W_m)=6g-3$. Since $H^3(W_m;\mathbb R)=0$, the map $H_1(Y_{g,4};\mathbb R)\to H_1(W_m;\mathbb R)$ is onto. Its kernel has dimension $6g-3$. The pair sequence identifies this kernel with the image of $H_2(W_m,Y_{g,4};\mathbb R)$, which also has dimension $6g-3$ by duality. Hence the preceding map $H_2(W_m;\mathbb R)\to H_2(W_m,Y_{g,4};\mathbb R)$ is zero. The real intersection pairing therefore vanishes. Since $H_2(W_m;\Z)$ is free, the integral intersection pairing vanishes as well.

\end{proof}

\begin{proposition}\label{prop:cover}
For $m\geq1$, the manifold $W_m$ has a connected regular cover with deck group $\Z/4\times\Z/2$ and
\[
 H_1(\widetilde W_m;\Z)=\Z^{8g}\oplus(\Z/m)^2.
\]
\end{proposition}
\begin{proof}

Let $a=a_1$, and let
\[
 \mathcal S=\{b_1,a_2,b_2,\ldots,a_g,b_g\}.
\]
The base group is free on $\{a\}\cup\mathcal S$. For the cyclic character
$a\mapsto1\in\mathbb Z/4$, $s\mapsto0$ for $s\in\mathcal S$, the Schreier
transversal $\{1,a,a^2,a^3\}$ gives the free basis
\[
 \alpha=a^4,\qquad \alpha_{j,s}=a^jsa^{-j},\quad
 0\leq j\leq3,\quad s\in\mathcal S
\]
of the kernel. There are $1+4(2g-1)=8g-3$ generators. The corresponding
base cover has four boundary components: the base-boundary word has zero
cyclic character, so every boundary lift has covering degree one.
Multiplicativity of Euler characteristic then gives genus $4g-3$.
The based monodromy of $\alpha$ is $r^4=\id$; that of $\alpha_{j,s}$ is the identity
unless $s=a_2$, in which case it is
\[
 r^j\tau^m r^{-j}.
\]

The connected double cover of the fiber is defined by sending every peripheral generator $x_j$ to $1\in\Z/2$. This is well defined because $x_1x_2x_3x_4=1$. Both $r$ and $\tau$ preserve this character. Each boundary component lifts to one circle of covering degree two. Since the covering fiber has Euler characteristic $2\chi(P)=-4$ and four boundary components, it has genus one.

We fix a point above the north pole and use the unique lifts of $r$ and $\tau$ fixing that point. Under the covering inclusion, these based lifts induce the restrictions of $r_*$ and $\tau_*$ to the parity kernel in the fiber group. Thus they are precisely the subgroup automorphisms used in the Schreier rewriting below. The lift of $r$ has fourth power equal to the identity: its fourth power is a deck transformation fixing that point. These lifts define the fiber double cover of the pulled-back bundle. Equivalently, on $\pi_1(W_m)$ the map to $\Z/4\times\Z/2$ sends $a_1$ to $(1,0)$, each $x_j$ to $(0,1)$, and the other base generators to zero. Its kernel defines the connected regular eight-fold cover $\widetilde W_m$.

The fiber group is free on $x_1,x_2,x_3$. The complete rewriting calculation is given in Appendix~\ref{app:schreier}. Applying Schreier rewriting with transversal $\{1,x_1\}$ to its parity kernel gives the free basis
\[
 q_1=x_2x_1^{-1},\quad q_2=x_3x_1^{-1},\quad
 q_3=x_1^2,\quad q_4=x_1x_2,\quad q_5=x_1x_3.
\]
In the corresponding homology basis, write $R$ and $S$ for the actions of the based lifts of $r$ and $\tau$. Reading the images of these five words gives
\[
 R=\begin{pmatrix}
 -1&-2&1&1&1\\
 1&0&0&0&-1\\
 0&-1&0&0&0\\
 0&0&1&0&-1\\
 0&-1&0&1&0
 \end{pmatrix},\qquad
 S=I+v\ell,
\]
where
\[
 v=(1,0,1,-1,0)^T,\qquad
 \ell=(0,1,0,0,-1).
\]

Schreier rewriting gives the following images in homology. The middle column gives the columns of $R$.
\begin{center}
\begin{tabular}{ccc}
\toprule
Class & Lift of $r$ & Lift of $\tau$\\
\midrule
$q_1$ & $-q_1+q_2$ & $q_1$\\
$q_2$ & $-2q_1-q_3-q_5$ & $q_2+v$\\
$q_3$ & $q_1+q_4$ & $q_3$\\
$q_4$ & $q_1+q_5$ & $q_4$\\
$q_5$ & $q_1-q_2-q_4$ & $q_5-v$\\
\bottomrule
\end{tabular}
\end{center}

For example, $r_*(q_1)=x_3x_2^{-1}$ has homology $-q_1+q_2$, giving the first column of $R$.

We have $R^4=I$ and $\ell v=0$, so $(v\ell)^2=0$ and $S^m=I+mv\ell$. The image lattice generated by the four conjugate increments $R^j(v\ell)R^{-j}$ is
\[
 L=\Z v+\Z Rv\subset\Z^5,
 \qquad Rv=(-1,1,0,1,-1)^T.
\]
Indeed, $\ell$ is primitive, so the image of the $j$th increment is $\Z R^jv$, and $R^2v=-v$, $R^3v=-Rv$. The first two rows of the matrix $(v,Rv)$ have determinant one. Thus $L$ is a primitive rank-two sublattice.

The coinvariants of the fiber-cover homology under the covered bundle monodromy are consequently
\[
 \Z^5/(mL)\cong\Z^3\oplus(\Z/m)^2.
\]
The north-pole section splits off the first homology of the base cover, which is $\Z^{8g-3}$. Therefore
\[
 H_1(\widetilde W_m;\Z)=\Z^{8g}\oplus(\Z/m)^2.
\]

\end{proof}

\begin{proof}[Proof of Theorem~\ref{thm:main}]
The filling, homology, support-genus, and covering assertions follow from Propositions~\ref{prop:fixedboundary}--\ref{prop:cover}, equation~\eqref{eq:boundaryH1}, and Corollary~\ref{cor:nonplanar}. It remains to distinguish infinitely many fundamental groups.

For a finitely generated abelian group $B$, its torsion order means the finite cardinality $|\Tor B|$, with value one when $\Tor B$ is trivial. A finitely generated group has finitely many subgroups of any fixed finite index, since each gives a permutation representation in a fixed finite symmetric group. If only finitely many isomorphism types occurred among the groups $\pi_1(W_m)$, the torsion orders of the abelianizations of all their index-eight subgroups would be bounded. The displayed covers have torsion order $m^2$, a contradiction. An infinite subfamily therefore has pairwise nonisomorphic fundamental groups, and its members are pairwise nonhomeomorphic.

\end{proof}

The argument does not require the chosen cover to be characteristic: a fixed fundamental group has only finitely many index-eight subgroups, and all of their abelianizations enter the bound.

The covering kernel gives a finer distinction when it is retained as part of the marking. Write $G_m=\pi_1(V_{g,m})$ and let
\[
 \chi_m:G_m\longrightarrow\Z/4\times\Z/2,
 \qquad K_m=\ker\chi_m,
\]
where $\chi_m(a_1)=(1,0)$, $\chi_m(x_j)=(0,1)$, and every other base generator maps to zero. These assignments respect the bundle relations and define an epimorphism for every $m\geq0$.

\begin{proposition}\label{prop:markedgroups}
Fix $g\geq2$ and $m,n\geq0$. If an isomorphism $G_m\to G_n$ carries $K_m$ onto $K_n$, then $m=n$.
\end{proposition}
\begin{proof}
Such an isomorphism identifies the abelianizations of the two kernels. For $m\geq1$, Proposition~\ref{prop:cover} gives
\[
 K_m/[K_m,K_m]\cong\Z^{8g}\oplus(\Z/m)^2.
\]
For $m=0$, the degree-four base pullback has trivial monodromy, and its fiber double cover is a product. The base has first Betti number $8g-3$, and the fiber double cover has first Betti number five. Hence
\[
 K_0/[K_0,K_0]\cong\Z^{8g+2}.
\]
The free ranks distinguish $m=0$ from every positive parameter. For positive $m,n$, equality of the torsion orders gives $m^2=n^2$, hence $m=n$.
\end{proof}

We can express a sufficient marking without naming the character. Let $F_m\lhd G_m$ be the fiber subgroup, let $S_m\leq G_m$ be the image of the north-pole section, and let $\mathcal P_m$ be the set of $F_m$-conjugacy classes of the four peripheral cyclic subgroups $\langle x_j\rangle$. An isomorphism preserving these data must send $F_m$ onto $F_n$, $\mathcal P_m$ onto $\mathcal P_n$, and $S_m$ onto a conjugate of $S_n$.

\begin{corollary}\label{cor:bundlemarked}
For fixed $g$, the groups $G_m$ with the fiber, peripheral, and section data just specified are pairwise nonisomorphic for distinct $m\geq0$.
\end{corollary}
\begin{proof}
The character $\eta_m:F_m\to\Z/2$ taking every peripheral generator to one is determined by $\mathcal P_m$: conjugating or inverting a peripheral generator does not change its value. It is invariant under the base monodromy, so it extends uniquely to a character $\epsilon_m:G_m\to\Z/2$ that vanishes on $S_m$.

Conjugation by $G_m$ permutes the four elements of $\mathcal P_m$. This permutation action has cyclic image of order four, generated by $a_1$; its kernel is the kernel of the first component of $\chi_m$. Therefore
\[
 K_m=\ker\bigl(G_m\to\operatorname{Sym}(\mathcal P_m)\bigr)
       \cap\ker\epsilon_m.
\]
An isomorphism preserving the stated data identifies the permutation kernels. After an inner conjugation in the target, it also sends $S_m$ onto $S_n$, and uniqueness identifies the two characters $\epsilon_m$ and $\epsilon_n$. Inner conjugation does not change either kernel. Proposition~\ref{prop:markedgroups} now gives $m=n$.
\end{proof}

This marking includes the fiber and section subgroups. We have not shown that the abstract group, or the boundary peripheral data alone, determines them.

\section{Fillings in arbitrary degree}\label{sec:general}
Fix $g\geq2$ and $p\geq4$. We use the two-vertex ribbon model of Section~\ref{sec:handles}, now with $p$ bands $e_0,\ldots,e_{p-1}$ in cyclic order. The rotation $r=r_p$ fixes the poles and cyclically permutes the bands. Set $z_j=e_je_0^{-1}$, with $z_0=z_p=1$. The based peripheral generators are
\[
 x_1=z_1^{-1},\qquad
 x_i=z_{i-1}z_i^{-1}\quad(2\leq i<p),\qquad
 x_p=z_{p-1}.
\]
Since $r_*(z_j)=z_{j+1}z_1^{-1}$, these satisfy
\[
 x_1\cdots x_p=1,\qquad r_*(x_j)=x_{j+1},
\]
with indices modulo $p$. In particular, $e_0e_2^{-1}=z_2^{-1}=x_1x_2$. As in the degree-four model, push this circle into the pair of pants containing the first two boundary components, leaving the north basepoint outside. This fixes the curve $\gamma$, its joining arc, and its based word $w=x_1x_2$. Choose its annular neighborhood disjoint from the pole and all boundary collars. The proof of Lemma~\ref{lem:standard-twist} gives a Dehn twist $\tau$ with
\[
 \tau_*(x_i)=wx_iw^{-1}\quad(i=1,2),\qquad
 \tau_*(x_i)=x_i\quad(3\leq i\leq p).
\]
In particular, $\tau_*(w)=w$.

Define $U_{g,p,m}\to B_g$ by the actual holonomies
\[
 a_1\longmapsto r,\qquad a_2\longmapsto\tau^{2m},
\]
with every other standard base generator mapping to the identity. We round its corners using the product collars. The boundary word has identity holonomy, and the collar representations are independent of $m$. Lemma~\ref{lem:boundary} gives Stein fillings of $(Y_{g,p},\xi_{g,p})$. Here $U_{g,p,0}=W_{g,p}$ and $U_{g,4,m}=V_{g,2m}$.

The pole gives a section, and both base and fiber are aspherical. Since $\tau$ acts trivially on $H_1(P_p;\Z)$, the cellular differential with fiber coefficients is the one in Proposition~\ref{prop:topology}. Thus
\[
 H_1(U_{g,p,m};\Z)=\Z^{2g}\oplus\Z/p,\quad
 H_2(U_{g,p,m};\Z)=\Z^{(2g-1)(p-1)},\quad H_k=0\ (k\geq3).
\]
The boundary has first Betti number $1+p(2g-1)$. Its difference from $b_1(U_{g,p,m})$ is $(2g-1)(p-1)=b_2(U_{g,p,m})$. The pair-sequence argument of Proposition~\ref{prop:ordinary} therefore makes the absolute intersection pairing zero. Corollary~\ref{cor:nonplanar}, applied to the cyclic filling $W_{g,p}$ of the same contact boundary, gives positive ordinary support genus.

Let $F=\pi_1(P_p)$ and set
\[
 A=H_1(P_p;\Z/2)\cong(\Z/2)^{p-1},\qquad
 F^{(2)}=\ker(F\to A),\qquad D=|A|=2^{p-1}.
\]
Write $\overline P\to P_p$ for this connected regular cover. Its first homology is free of rank
\[
 h_p=1+D(p-2).
\]
Indeed, a spine of $P_p$ is a rose with $p-1$ edges, and its cover has $D$ vertices and $D(p-1)$ edges. Both $r$ and $\tau$ preserve $F^{(2)}$. We always use their lifts fixing a chosen point above the north pole.

\begin{lemma}\label{lem:general-lift}
Let $T$ be the action of the based lift of $\tau^2$ on $H=H_1(\overline P;\Z)$. Then $T=I+N$, where $N^2=0$ and $N\neq0$. If $R$ is the action of the based lift of $r$, then $R^p=I$.
\end{lemma}
\begin{proof}
In $A$, the class of $\gamma$ is $\bar x_1+\bar x_2\neq0$, of order two. Every annulus above the twisting annulus consequently covers it with degree two. The based lift of $\tau^2$ is a product of Dehn twists, all of the same sign, about the disjoint lifts of $\gamma$. To see that it is this lift, write the downstairs twist as $(\theta,s)\mapsto(\theta+2\pi h(s),s)$. On a covering annulus with $\theta=2\widetilde\theta$, the square lifts to $(\widetilde\theta,s)\mapsto(\widetilde\theta+2\pi h(s),s)$ and extends by the identity off these annuli. It fixes the chosen point. For each lifted twist curve $\delta_i$, let $N_i$ be the signed homological transvection increment, so $N_i(a)=\pm(a\cdot[\delta_i])[\delta_i]$. The lifted curves are disjoint, and their pairwise algebraic intersections, including their self-intersections, vanish. Hence $N_iN_j=0$ for all $i,j$. The action of their product of twists is therefore $T=\prod_i(I+N_i)=I+N$, where $N=\sum_iN_i$ satisfies $N^2=0$.

We verify $N\neq0$ on the covering graph. Use the free generators $x_1,\ldots,x_{p-1}$ of $F$ and the corresponding basis $u_1,\ldots,u_{p-1}$ of $A$. Its vertices are $a\in A$, with an oriented edge $E(a,i)$ from $a$ to $a+u_i$ for each $i$. Distinct lifts of the same rose edge are distinct edges, even when they have the same endpoints.

For $a\in A$, let $P(a)$ be the chain of the closed path $w^2$ starting at $a$. Explicitly,
\[
 P(a)=E(a,1)+E(a+u_1,2)
       +E(a+u_1+u_2,1)+E(a+u_2,2).
\]
The word $c=(x_1x_3)^2$ also lifts to a closed path at zero. Since
\[
 \tau_*^2(x_i)=w^2x_iw^{-2}\ (i=1,2),\qquad
 \tau_*^2(x_3)=x_3,
\]
the difference between the chains of its image and itself is
\[
 P(0)-P(u_1)+P(u_1+u_3)-P(u_3).
\]
The coefficient of $E(0,2)$ is $-1$. Here $u_1,u_2,u_3$ are independent because $p\geq4$. This is a nonzero cycle in a graph, which has no cellular two-boundaries. It represents $N[c]\neq0$.

Finally, the $p$th power of the based lift of $r$ is a deck transformation fixing the chosen point, so it is the identity. Hence $R^p=I$.
\end{proof}

The first homology of a fixed-degree cover now detects the variation. Let $\bar r$ be the action of $r$ on $A$, and form
\[
 Q_p=A\rtimes_{\bar r}\Z/p.
\]
There is an epimorphism
\[
 \Psi_m:\pi_1(U_{g,p,m})\longrightarrow Q_p,
\]
which sends $x_i$ to $(\bar x_i,0)$, $a_1$ to $(0,1)$, and the other base generators to zero. The fiber relation holds in $A$, the rotation relations give its stated action, and the twist relations hold because $\tau$ acts trivially on $A$. Thus the associated cover $\widetilde U_{g,p,m}$ is connected and regular of degree $pD$.

The base of this cover is the degree-$p$ cyclic cover of $B_g$ defined by $a_1\mapsto1$, with every other base generator mapping to zero. Its fundamental group is free of rank
\[
 b_{g,p}=1+p(2g-1).
\]
The Schreier generators are $a_1^p$ and $a_1^jsa_1^{-j}$, where $0\leq j<p$ and $s\in\{b_1,a_2,b_2,\ldots,a_g,b_g\}$. On the covering fiber, their based homological monodromies are the identity except for $s=a_2$, when they are
\[
 R^jT^mR^{-j}=I+mR^jNR^{-j}.
\]
The north-pole section lifts over this base cover and splits off its first homology. Consequently,
\[
 H_1(\widetilde U_{g,p,m};\Z)
 \cong\Z^{b_{g,p}}\oplus H/(mL_p),\qquad
 L_p=\sum_{j=0}^{p-1}\operatorname{im}(R^jNR^{-j})\subset H.
\]
Each conjugate increment $R^jNR^{-j}$ is independent of $m$. Since $T^m=I+mN$, the image of the corresponding monodromy minus the identity is $m\operatorname{im}(R^jNR^{-j})$. Summing these relation subgroups gives exactly $mL_p$. Abelianizing the split extension of the free base group by $F^{(2)}$ then gives the displayed coinvariant quotient.

Let $s_p=\operatorname{rank}L_p$. Lemma~\ref{lem:general-lift} gives $s_p\geq1$. The integral Smith normal form of the inclusion $L_p\subset H$ has positive factors $d_1,\ldots,d_{s_p}$, depending only on $p$. No primitivity of $L_p$ is needed. For $m\geq1$ we obtain
\begin{equation}\label{eq:general-cover}
 H_1(\widetilde U_{g,p,m};\Z)
 \cong\Z^{b_{g,p}+h_p-s_p}
       \oplus\bigoplus_{i=1}^{s_p}\Z/(m d_i).
\end{equation}
In particular, its torsion order is $c_p m^{s_p}$, where $c_p=\prod_i d_i>0$. For $m=0$ the cover is a product over the cyclic base cover, and its first homology is $\Z^{b_{g,p}+h_p}$.

\begin{proof}[Proof of Theorem~\ref{thm:general}]
The preceding construction proves the filling and ordinary-topology assertions. Equation~\eqref{eq:general-cover} supplies covers of the fixed degree $p2^{p-1}$ whose torsion orders are unbounded. The finite-index subgroup argument in the proof of Theorem~\ref{thm:main} therefore gives infinitely many isomorphism types of fundamental groups.

For the ruled realization, choose an $r$-invariant sphere area form and represent $\tau$ by an area-preserving annular twist. Both maps fix the pole, and $\tau$ is the identity near it and near the capping disks. Capping the boundary components gives a flat sphere bundle over $B_g$. The actual holonomies satisfy $\prod_{i=1}^g[\rho(a_i),\rho(b_i)]=\id$ on the capped sphere, so the representation factors through $\pi_1(\Sigma_g)$. Its flat sphere bundle over the closed base restricts to the capped bundle over $B_g$, giving the required flat extension across the missing base disk. The pole section has flat unitary normal holonomy, a rotation along $a_1$ and the identity along the other generators. Its normal Euler class is zero, since $H^2(\Sigma_g;\Z)$ is torsion-free. The section-preserving clutching argument in Proposition~\ref{prop:ruled} identifies the sphere bundle with $\Sigma_g\times S^2$.

The centers of the capping disks give a connected unbranched multisection, since $r$ permutes them transitively. The sum of the invariant fiber area form and a positive base area form is symplectic and makes this multisection symplectic. Its covering degree gives genus $1+p(g-1)$. Its disjointness from the pole section fixes its class as $p[\Sigma_g\times\{\mathrm{pt}\}]$, of square zero. Removing its neighborhood and one sphere fiber recovers $U_{g,p,m}$.
\end{proof}

Equation~\eqref{eq:general-cover} also distinguishes every parameter when the kernel of $\Psi_m$ is retained. An isomorphism preserving that kernel preserves its abelianization. The free rank separates $m=0$ from $m>0$, and the torsion order $c_p m^{s_p}$ separates the positive parameters. We do not assert that every abstract isomorphism preserves this kernel.

The restriction $p\geq4$ is needed for this variation. When $p=3$, the relation $x_1x_2x_3=1$ makes the same twist act by conjugation by $w=x_3^{-1}$ on the entire fiber group. Its powers can be absorbed into the second base generator, so the resulting abstract groups are isomorphic. When $p=2$, the fiber group is cyclic and the annular twist acts trivially on it. The cyclic complements in these two degrees remain covered by Proposition~\ref{prop:topology}, but this twist does not produce infinitely many isomorphism types of fundamental groups in these degrees.

\section{Ruled realization and further questions}\label{sec:ruled}
\begin{proposition}\label{prop:ruled}
For every $g\geq2$ and $m\geq0$, $V_{g,m}$ is diffeomorphic to the complement of a connected symplectic degree-four unbranched multisection and a regular fiber in a symplectic ruled surface diffeomorphic to $\Sigma_g\times S^2$. The multisection has genus $4g-3$, square zero, and homology class $4[\Sigma_g\times\{\mathrm{pt}\}]$.
\end{proposition}

\begin{proof}
We choose the actual maps $r$ and $\tau$ to preserve a common area form $\sigma$ on $S^2$. Take $\sigma$ invariant under the standard rotation $r$. In area coordinates $(\theta,s)$ on an annular neighborhood of the standard curve $\gamma$, with $\theta\in\mathbb R/2\pi\Z$ and $\sigma=ds\wedge d\theta$, a Dehn twist has the form
\[
 (\theta,s)\longmapsto(\theta+2\pi h(s),s),
\]
where $h$ is constant zero near one end and constant one near the other. This map preserves $\sigma$ and is the identity near both ends. We extend it by the identity outside the annulus and choose its sign to agree with the based action in Section~\ref{sec:variation}. Its support misses the disk centers, the north pole, and the boundary collars.

Cap the four boundary components by disks. The maps $r$ and $\tau$ extend to the capped sphere, and their boundary-word holonomy is the identity. The flat sphere bundle therefore extends across the missing base disk to an oriented $S^2$-bundle over $\Sigma_g$. The north pole gives a section. Its normal holonomy is a rotation along $a_1$ and the identity on the other generators. Hence its normal oriented plane bundle has a flat unitary connection. Its Euler class has zero image in real cohomology and therefore vanishes in $H^2(\Sigma_g;\Z)$.

By Smale's theorem, $SO(3)\hookrightarrow\operatorname{Diff}^+(S^2)$ is a homotopy equivalence \cite{Smale}, so the sphere bundle admits an $SO(3)$ reduction. The pole section gives a unit vector in the associated oriented rank-three bundle. Its orthogonal complement is the normal plane bundle of the section. The vanishing Euler class trivializes this plane bundle; a frame for it and the unit vector then trivialize the rank-three bundle. The induced sphere-bundle trivialization takes the pole section to $\Sigma_g\times\{\mathrm{pt}\}$ in $\Sigma_g\times S^2$.

The four disk centers form a connected unbranched degree-four multisection $F_m$, since $r$ permutes them transitively and $\tau$ fixes them. The covering formula gives $g(F_m)=4g-3$. Write its class as $4[\Sigma_g\times\{\mathrm{pt}\}]+k[\{\mathrm{pt}\}\times S^2]$. Its disjointness from the pole section gives $k=0$. Thus
\[
 [F_m]=4[\Sigma_g\times\{\mathrm{pt}\}],\qquad F_m^2=0.
\]

In the flat model, the forms $\sigma$ on the sphere and a positive area form on the base descend to closed forms. Their sum is symplectic and restricts positively to $F_m$, since the disk centers are locally constant in flat charts. Removing a neighborhood of $F_m$ and one regular sphere fiber recovers $V_{g,m}$.
\end{proof}

The pure-braid variation preserves the marked spinal boundary because its support misses the fiber collars and its boundary-word holonomy is the identity. The closed-sum maps in \cite{AY} are specified on a different boundary. To use those maps for fillings, one must first identify the corresponding spinal markings and verify that the contact structure is unchanged.

\begin{question}
Can the closed-sum gluing maps of \cite{AY} be realized by variations of Stein fillings of a fixed $(Y_{g,4},\xi_{g,4})$?
\end{question}

Proposition~\ref{prop:markedgroups} recovers the parameter from the group with its covering kernel. Corollary~\ref{cor:bundlemarked} gives a sufficient fiber, peripheral, and section marking. The unmarked classification remains open.

Here the boundary peripheral data mean the inclusion-induced homomorphism
\[
 \iota_m:\pi_1(Y_{g,4})\longrightarrow G_m,
\]
considered up to conjugation in $G_m$ and precomposition by automorphisms induced by orientation-preserving self-diffeomorphisms of $Y_{g,4}$. An isomorphism $\varphi:G_m\to G_n$ preserves these data if $\varphi\circ\iota_m$ agrees, up to conjugation in $G_n$, with $\iota_n\circ h_*$ for such a boundary diffeomorphism $h$. These data retain the full boundary homomorphism and are distinct from the fiber peripheral classes $\mathcal P_m$; they do not include the fiber and section subgroups as additional markings.
\begin{question}
For a fixed $g$, when are $G_m$ and $G_n$ isomorphic? Does every isomorphism preserving the boundary peripheral data carry $K_m$ onto $K_n$?
\end{question}

The existence of exotic Stein fillings with prescribed fundamental group is known. Akhmedov and Ozbagci \cite[Theorem 1]{AkO} prove that every finitely presented group occurs as the fundamental group of infinitely many pairwise homeomorphic but nondiffeomorphic Stein fillings of the canonical contact structure on a link of an isolated complex surface singularity. Applied to $G_m$, their theorem gives such a contact manifold for each fundamental group in our family. It does not identify that contact manifold with $(Y_{g,4},\xi_{g,4})$, or identify the homeomorphism type or asphericity of those fillings with $V_{g,m}$. Akbulut and Yasui \cite{AkY} give simply connected examples with second Betti number two. Yasui's partial-twist construction \cite{Yasui} also requires a prescribed configuration of vanishing cycles in a positive allowable Lefschetz fibration over the disk; its algorithm may change the original boundary. Thus these results do not settle the following question.

\begin{question}
Does $(Y_{g,4},\xi_{g,4})$ admit homeomorphic but nondiffeomorphic Stein fillings?
\end{question}

The zero intersection pairing restricts the local constructions available inside our fillings.

\begin{proposition}\label{prop:subdomainobstruction}
Let $X$ be a compact oriented smooth four-manifold embedded as a codimension-zero submanifold of $V_{g,m}$. Then the absolute intersection pairing of $X$ is zero. Moreover, every map $S^2\to V_{g,m}$ represents zero in second homology.
\end{proposition}
\begin{proof}
Represent classes in $H_2(X;\Z)$ by cycles supported away from its boundary and perturb them to be transverse. Their intersections computed in $X$ and in $V_{g,m}$ agree. For the inclusion $i$, this gives
\[
 Q_X(u,v)=Q_{V_{g,m}}(i_*u,i_*v)=0.
\]
The second assertion follows from asphericity: every sphere map is nullhomotopic.
\end{proof}

In particular, a nucleus with nonzero intersection pairing cannot embed as a codimension-zero submanifold of $V_{g,m}$. Neither can a neighborhood of two closed oriented surfaces with nonzero algebraic intersection, or of a surface with nonzero self-intersection. This does not exclude an embedded contractible cork, whose second homology vanishes.

For background on cork decompositions, see \cite[Chapter~10]{AkBook}. The following capping criterion distinguishes relative from absolute nondiffeomorphism. Here diffeomorphisms and homeomorphisms preserve orientation.

\begin{proposition}\label{prop:corkcriterion}
Let $V$ be a Stein filling of $(Y,\xi)$, let $C\subset\Int V$ be a compact contractible smooth submanifold, and let $f:\partial C\to\partial C$ be a diffeomorphism extending to a homeomorphism of $C$. Form
\[
 V_f=(V\setminus\Int C)\cup_f C
\]
with the boundary identified with $Y$ by the unchanged exterior. Suppose $V_f$ admits a Stein structure inducing $\xi$ under this boundary identification. Let $Z$ be a compact smooth four-manifold with boundary $-Y$. If the closed manifolds $V\cup_Y Z$ and $V_f\cup_Y Z$ are not diffeomorphic, then $V$ and $V_f$ are homeomorphic relative to $Y$ but are not diffeomorphic relative to $Y$. If every orientation-preserving self-diffeomorphism of $Y$ extends over $Z$, then $V$ and $V_f$ are not orientation-preservingly diffeomorphic even without fixing the boundary. In addition, if $V$ is aspherical, then $V_f$ is aspherical.
\end{proposition}
\begin{proof}
A topological extension of $f$ on $C$, together with the identity on the exterior, gives a homeomorphism $V_f\to V$ fixed on $Y$. A diffeomorphism relative to $Y$ would extend by the identity on $Z$, contradicting the hypothesis on the closed manifolds. More generally, the restriction of any orientation-preserving diffeomorphism $V\to V_f$ is a self-diffeomorphism of $Y$. If it extends over $Z$, the two maps glue to a diffeomorphism of the closed manifolds, again a contradiction. Finally, asphericity is preserved by the homeomorphism.
\end{proof}

To apply Proposition~\ref{prop:corkcriterion} to $V_{g,m}$, one still needs a cork embedding, a Stein structure after twisting with the prescribed boundary contact structure, and a smooth invariant distinguishing the capped manifolds. Absolute nondiffeomorphism also requires control of the boundary diffeomorphisms, unless it can be detected directly on the fillings. The handle construction in Section~\ref{sec:handles} does not yet supply the marked Legendrian front for this calculation.

Corollary~\ref{cor:nonplanar} gives a lower bound of one for the ordinary support genus. We do not know its exact value.
\begin{question}
What is the ordinary support genus of $\xi_{g,4}$?
\end{question}

\appendix
\section{Schreier rewriting in the fiber double cover}\label{app:schreier}
We compute the actions of the based lifts on the fiber-cover group, whose abelianizations give the matrices in Proposition~\ref{prop:cover}.
Write $x=x_1$, $y=x_2$, and $z=x_3$. The fiber group is free on $x,y,z$;
the fourth peripheral generator is $(xyz)^{-1}$. The two automorphisms are
\[
\begin{aligned}
 r(x)&=y,& r(y)&=z,& r(z)&=z^{-1}y^{-1}x^{-1},\\
 \tau(x)&=xyxy^{-1}x^{-1},&
 \tau(y)&=xyx^{-1},& \tau(z)&=z.
\end{aligned}
\]
Here the twist formulas follow from Lemma~\ref{lem:standard-twist}, with
$w=xy$. In particular, both automorphisms preserve the character
$\eta:\mathbb F(x,y,z)\to\mathbb Z/2$ sending $x,y,z$ to one.

Use the right-coset transversal $\{1,x\}$ for $\ker\eta$. The positive
Schreier edges are
\[
\begin{array}{c|ccc}
\text{initial parity}&x&y&z\\ \hline
0&1&q_1&q_2\\
1&q_3&q_4&q_5
\end{array}
\qquad
\begin{aligned}
 q_1&=yx^{-1},&q_2&=zx^{-1},\\
 q_3&=x^2,&q_4&=xy,&q_5&=xz.
\end{aligned}
\]
Each letter changes parity. To read a negative letter $t^{-1}$ from parity
$\epsilon$, use the inverse of the positive Schreier edge for $t$ starting
at parity $1-\epsilon$, then change parity. Start at parity zero. An even
word ends at parity zero and rewrites as a word in $q_1,\ldots,q_5$.
Indeed the successive factors are
$\overline{u}\,t\,\overline{ut}^{-1}$; their transversal terms telescope.
The omitted edge is the tree edge labelled $x$ from parity zero to parity
one. The other five edges give the displayed free basis of the two-sheeted
covering graph.

Applying this rule to the images of the five generators gives exact
identities in the parity kernel:
\[
\begin{array}{c|l|l}
 & r_*(q_i)&\tau_*(q_i)\\ \hline
q_1&q_2q_1^{-1}&q_4q_1q_4^{-1}\\
q_2&q_5^{-1}q_1^{-1}q_3^{-1}q_1^{-1}&q_2q_3q_1q_4^{-1}\\
q_3&q_1q_4&q_4q_3q_4^{-1}\\
q_4&q_1q_5&q_4\\
q_5&q_1q_2^{-1}q_4^{-1}&q_4q_1^{-1}q_3^{-1}q_5
\end{array}
\]
For example,
\[
 r_*(q_2)=r(zx^{-1})=z^{-1}y^{-1}x^{-1}y^{-1}.
\]
Its four negative letters, read from parity zero, contribute respectively
$q_5^{-1},q_1^{-1},q_3^{-1},q_1^{-1}$. For the twist,
\[
 \tau_*(q_2)=zxyx^{-1}y^{-1}x^{-1}.
\]
The successive contributions are $q_2,q_3,q_1,1,q_4^{-1},1$, yielding
$q_2q_3q_1q_4^{-1}$. Each other entry can also be checked by substituting
the defining words for the $q_i$ and freely reducing.

The table uses the lifts fixing the chosen point above the north pole. Its abelianization gives the matrices $R$ and $S$ in
Proposition~\ref{prop:cover}. In particular,
\[
 S-I=v\ell,\quad
 v=(1,0,1,-1,0)^T,\quad\ell=(0,1,0,0,-1),
 \qquad R^4=I,
\]
and direct multiplication gives
\[
 \ell v=0,\qquad Rv=(-1,1,0,1,-1)^T,\qquad R^2v=-v.
\]
Thus the four lifted monodromies have increments
$mR^j(v\ell)R^{-j}$. Their images generate
$m(\mathbb Zv+\mathbb ZRv)$. The coefficient row $\ell R^{-j}$ is primitive
because $R$ is unimodular, and the first two rows of $(v,Rv)$ have
determinant one. This verifies the integral quotient
\[
 \mathbb Z^5\big/m(\mathbb Zv+\mathbb ZRv)
 \cong\mathbb Z^3\oplus(\mathbb Z/m)^2
\]
used in the covering calculation.

\section{A direct holomorphic quotient construction}\label{app:quotient}
The invariant functions $f_B,f_P$ of Proposition~\ref{prop:stein} give a strictly pseudoconvex domain in $E$ diffeomorphic to $W_{g,p}$. The supported boundary contact structure $\xi_{g,p}$ is defined using the additive potential in Proposition~\ref{prop:stein}.

On the product consider the invariant function
\[
 \Phi_\epsilon=\epsilon\log\left(e^{f_B/\epsilon}+e^{f_P/\epsilon}\right),
 \qquad\epsilon>0.
\]
It is strictly plurisubharmonic. Indeed, if
\[
 \theta=\frac{e^{f_B/\epsilon}}{e^{f_B/\epsilon}+e^{f_P/\epsilon}},
\]
its Levi form is the sum of $\theta\mathcal L(f_B)+(1-\theta)\mathcal L(f_P)$ and a positive semidefinite form. Both weights are positive. The two Levi forms are positive in their respective factor directions, so this sum is positive on every nonzero complex tangent vector.

Choose $\delta>0$ so small that $-2\delta\leq f_B\leq0$ and $-2\delta\leq f_P\leq0$ are boundary collars with nonvanishing outward derivative. Take $\epsilon\log2<\delta/2$ and $c=-\delta$. The inequalities
\[
 \max(f_B,f_P)\leq\Phi_\epsilon
 \leq\max(f_B,f_P)+\epsilon\log2
\]
show that on $\Phi_\epsilon=c$, at least one factor lies in its boundary collar. Choose invariant outward collar vector fields on the two factors, cut off outside these collars, with nonnegative derivative on their factor functions. Their sum has strictly positive derivative on $\Phi_\epsilon=c$, since both weights in $d\Phi_\epsilon=\theta\,df_B+(1-\theta)\,df_P$ are positive. Thus the level is regular. Moving along the collars identifies its sublevel equivariantly with an inward truncation and smoothing of the product corners. Taking the quotient gives a domain diffeomorphic to $W_{g,p}$.

The function descends to the quotient. With its complex structure $J$, define
\[
 \lambda=-d\Phi_\epsilon\circ J,\qquad\omega=d\lambda.
\]
Strict plurisubharmonicity makes $\omega$ symplectic. The Liouville vector field is the gradient of $\Phi_\epsilon$ for the associated K\"ahler metric and points outward on the sublevel boundary. Its boundary restriction is therefore a positive contact form. The domain is Stein; its interior has the proper strictly plurisubharmonic exhaustion $-\log(c-\Phi_\epsilon)$.

\section*{Acknowledgments}
The author thanks B. Doug Park, Sai-Kee Yeung, and \c{C}a\u{g}r{\i} Karakurt for helpful discussions on braided multisections. He is grateful to B. Doug Park for several discussions between 2007 and 2010 on ideas for constructing fake quadrics via braided multisections. He also thanks Weiyi Zhang for their collaboration in 2011 on braided multisections and the summing-spheres construction, and Burak Ozbagci for their collaboration in 2012 on exotic Stein fillings \cite{AkO}. The author discussed the braided-multisection construction in a talk at the Istanbul Center for Mathematical Sciences (IMBM) in July 2015. He thanks \c{C}a\u{g}r{\i} Karakurt and Ferit \"{O}zt\"{u}rk for the invitation to IMBM.

An LLM-based tool was used for assistance with computational preparation, grammar checking, and figure generation. All mathematical ideas and results are due to the author, who takes full responsibility for the content.

\medskip
\noindent\begin{minipage}{\linewidth}\small
\textsc{Anar Akhmedov}\\
School of Mathematics, University of Minnesota\\
Visiting Scholar, Department of Mathematics, Harvard University\\
\href{mailto:akhmedov@umn.edu}{\texttt{akhmedov@umn.edu}}
\end{minipage}

\begin{thebibliography}{99}\small
\setlength{\itemsep}{1pt}
\bibitem{AkBook} S. Akbulut, \emph{$4$-Manifolds}, Oxford University Press, Oxford, 2016.
\bibitem{AO} S. Akbulut and B. Ozbagci, \emph{On the topology of compact Stein surfaces}, Int. Math. Res. Not. 2002, no. 15, 769--782; \href{https://arxiv.org/abs/math/0103106}{arXiv:math/0103106}.
\bibitem{AkY} S. Akbulut and K. Yasui, \emph{Infinitely many small exotic Stein fillings}, J. Symplectic Geom. 12 (2014), no. 4, 673--684.
\bibitem{AB} A. Akhmedov, \emph{Braided multisections and symplectic four-manifolds with the rational cohomology of $S^2\times S^2$}, manuscript, 2026, \href{https://arxiv.org/abs/2609.11727}{arXiv:2609.11727}.
\bibitem{I} A. Akhmedov, \emph{Topology of Akhmedov's Fake Quadrics}, manuscript.
\bibitem{II} A. Akhmedov, S. Akbulut, and E. Z. Y{\i}ld{\i}z, \emph{Handlebodies for Akhmedov's Fake Quadrics}, manuscript.
\bibitem{AEMS} A. Akhmedov, J. B. Etnyre, T. E. Mark, and I. Smith, \emph{A note on Stein fillings of contact manifolds}, Math. Res. Lett. 15 (2008), no. 6, 1127--1132; \href{https://arxiv.org/abs/0712.3932}{arXiv:0712.3932}.
\bibitem{AkO} A. Akhmedov and B. Ozbagci, \emph{Exotic Stein fillings with arbitrary fundamental group}, Geom. Dedicata 195 (2018), 265--281; \href{https://arxiv.org/abs/1212.1743}{arXiv:1212.1743}.
\bibitem{APOdd} A. Akhmedov and B. D. Park, \emph{Exotic smooth structures on small $4$-manifolds with odd signatures}, Invent. Math. 181 (2010), no. 3, 577--603; \href{https://arxiv.org/abs/math/0701829}{arXiv:math/0701829}.
\bibitem{AY} A. Akhmedov and S.-K. Yeung, \emph{Braided Multisection Sums, Cyclic Covers, and Fake Quadrics}, manuscript, 2026.
\bibitem{AZ} A. Akhmedov and W. Zhang, \emph{The fundamental group of symplectic $4$-manifolds with $b^+=1$}, preprint, 2011; posted on arXiv in 2015, \href{https://arxiv.org/abs/1506.08367}{arXiv:1506.08367}.
\bibitem{BVHM} R. I. Baykur and J. Van Horn-Morris, \emph{Fillings of genus-1 open books and 4-braids}, Int. Math. Res. Not. IMRN 2018, no. 5, 1329--1346; \href{https://arxiv.org/abs/1604.02945}{arXiv:1604.02945}.
\bibitem{CE} K. Cieliebak and Y. Eliashberg, \emph{From Stein to Weinstein and Back: Symplectic Geometry of Affine Complex Manifolds}, American Mathematical Society Colloquium Publications, vol. 59, American Mathematical Society, Providence, RI, 2012.
\bibitem{Etnyre} J. B. Etnyre, \emph{Planar open book decompositions and contact structures}, Int. Math. Res. Not. 2004, no. 79, 4255--4267; \href{https://arxiv.org/abs/math/0404267}{arXiv:math/0404267}.
\bibitem{Gompf} R. E. Gompf, \emph{Handlebody construction of Stein surfaces}, Ann. of Math. (2) 148 (1998), 619--693.
\bibitem{HP} C. S. Hays and B. D. Park,
\emph{Non-isotopic symplectic surfaces in product 4-manifolds},
Trans. Amer. Math. Soc. 360 (2008), no. 11, 5771--5788;
\href{https://doi.org/10.1090/S0002-9947-08-04717-X}{doi:10.1090/S0002-9947-08-04717-X}.
\needspace{18\baselineskip}
\bibitem{LVW} S. Lisi, J. Van Horn-Morris, and C. Wendl, \emph{On symplectic fillings of spinal open book decompositions I: Geometric constructions}, \href{https://arxiv.org/abs/1810.12017}{arXiv:1810.12017v2}, 2026.
\bibitem{Smale} S. Smale, \emph{Diffeomorphisms of the $2$-sphere}, Proc. Amer. Math. Soc. 10 (1959), no. 4, 621--626.
\bibitem{Smith} I. Smith, \emph{Symplectic submanifolds from surface fibrations}, Pacific J. Math. 198 (2001), no. 1, 197--205.
\bibitem{Yasui} K. Yasui, \emph{Partial twists and exotic Stein fillings}, preprint, 2014, \href{https://arxiv.org/abs/1406.0050}{arXiv:1406.0050}.
\bibitem{Zhou} Z. Zhou, \emph{A note on contact manifolds with infinite fillings}, \href{https://arxiv.org/abs/2304.11804}{arXiv:2304.11804}, 2023.
\end{thebibliography}
\end{document}